\documentclass[12pt,a4paper,twoside,final,notitlepage, reqno]{article}

\usepackage[english]{babel}
\usepackage[latin1]{inputenc}
\usepackage[a4paper,left=3.5cm,right=3.5cm]{geometry}
\usepackage{amssymb}
\usepackage{amsthm} 
\usepackage{amsmath}
\usepackage{mathtools}
\usepackage{amscd}
\usepackage{geometry}
\usepackage{graphics,graphicx}
\usepackage{epstopdf}
\usepackage[usenames, dvipsnames]{color}
\usepackage[textsize=small]{todonotes}
\usepackage{booktabs}
\usepackage{subfigure} 
\usepackage{siunitx}

\graphicspath{{figs/}}
\makeatletter
\def\input@path{{figs/}}
\makeatother

\theoremstyle{definition}
\newtheorem{theorem}{Theorem}%[section]

\newtheorem{proposition}{Proposition}
\newtheorem{corollary}{Corollary}

\theoremstyle{definition}
\newtheorem{definition}{Definition}%[section]

\theoremstyle{definition}
\newtheorem{remark}{Remark}%[section]

\providecommand{\keywords}[1]  {\textbf{Keywords:} #1}
\providecommand{\subjclass}[1] {\textbf{Subject classification:} #1}

\let\div\undefined
\DeclareMathOperator{\div}{div}

\usepackage{esdiff}
\renewcommand{\d}[1]{\mathrm{d} #1}
\def \dx {\d{x}}

\def \dl {\d{l}}

\newcommand{\N}{\mathbb{N}}
\newcommand{\Z}{\mathbb{Z}}

\newcommand{\R}{\mathbb{R}}

\renewcommand{\phi}{\varphi}

\usepackage[all]{xy}

\def\di{\displaystyle}

\newcommand \meas[1]   {\vert {#1} \vert }

\def \msh    {\mathcal{T}}   % maillage
\def \fc     {\mathcal{E}}   % faces (ou aretes)
\def \bfc    {\fc_0}         % boundary faces 
\def \ifc    {\fc_i}         % internal faces
\def \vtc    {\mathcal{N}}   % vertexes 

\def \LL   {\mathcal{L}}   % Fonctionelle lagrangienne
\def \HH   {\mathcal{H}}   % Fonctionelle hamiltonieinne

\def\div{\mathop{\rm div}\nolimits}

\def\bj {{\bf j}}
\def\bn {{\bf n}}
\def \pp    {\mathbf{p}}
\def \qq    {\mathbf{q}}
\usepackage{authblk}

\title{\bf{\Large{
    Discrete embeddings 
    \\[5pt]
    for Lagrangian  and Hamiltonian systems}}
}

\author[1]{Jacky Cresson \thanks{jacky.cresson@univ-pau.fr}}
\author[1]{Isabelle Greff \thanks{isabelle.greff@univ-pau.fr}}
\author[1]{Charles Pierre \thanks{charles.pierre@univ-pau.fr}}

\affil[1]{
  Laboratoire  de Math\'ematiques et de leurs Applications, 
  UMR CNRS 5142, \protect\\
  Universit\'e de Pau et des Pays de l'Adour, France}

\usepackage{fancyhdr}

\begin{document}

\date{16 January, 2018}

\maketitle

\noindent
\\ \\ 
 \\ \\
% \acknow{}
%% 
%% 
%% 
%% 
%% 
\begin{abstract}
  The  topic of this paper is to study the conservation of variational properties for a given problem  when discretising it.
  Precisely 
  we are interested in  Lagrangian or Hamiltonian structures and thus with variational problems attached to a least action principle.
  Consider a partial
  differential equation (PDE) deriving from  a variational 
  principle. 
  A natural question
  is to know whether this structure 
  is preserved at the discrete level when discretising the PDE.
  To address this question a concept of \textit{coherence} is introduced. 
  Both the differential equation (the PDE translating the least action principle) and the variational structure can be embedded at the discrete level. 
  This provides two discrete embeddings for the original problem. 
  If these procedures finally provide the same discrete problem we will say that the discretisation is \textit{coherent}.
  Our purpose is illustrated with the Poisson problem.
  Coherence 
  for discrete embeddings of Lagrangian structures
  is studied for various classical discretisations.
  For Hamiltonian structures, we show the coherence between a discrete Hamiltonian and the discretisation of the mixed formulation of the Poisson problem.
\end{abstract}
%% 
%% 
%% 
%% 

%%%%%%%%%%%%%%%%%%%%%%%%%%%% 
%%%%%%%%%%%%%%%%%%%%%%%%%%%% 
\section*{Introduction}
%%%%%%%%%%%%%%%%%%%%%%%%%%%% 
%%%%%%%%%%%%%%%%%%%%%%%%%%%% 
%% 
%% 
% \label{}
%% 
%% 
%% 
%% 
%% 
%% 
%% 
%% 
%% 
%% 
%% 
%% 
Many problems in physics, formulated in terms of Partial Differential Equations (PDE), are associated with essential  structural properties.
For instance we mention the maximum principle, conservation laws or variational principles in mechanics.
It is quite natural to ask the numerical methods to preserve these structural properties at the discrete level:  in order to enforce the numerical solutions to satisfy the underlying physics of the problem.
\\ \\
Two fundamental notions arising in classical mechanics are Lagrangian and Hamiltonian structures.
Lagrangian systems are made of one functional, called the Lagrangian functional, and a variational principle called the least action principle. 
From the least action principle is derived a second order 
differential equation called the Euler-Lagrange equation, see e.g. \cite{Ar}.
The Lagrangian structure is much more fundamental than its associated Euler-Lagrange equation: it contains information that the Euler-Lagrange equation does not. 
An important example is the change of coordinates. 
The Lagrangian structure is independent from change of coordinates, whereas the associated Euler-Lagrange equation may completely change of nature 
(from linear to non linear for instance).
% \\
Similarly, Hamiltonian systems also are associated to a
variational structure. % (precisely described in section \ref{part3}).
They
are associated with fundamental properties such as energy conservation or existence of first integrals.
\\ \\
Consider a numerical method for the resolution of  a problem that derives from a variational principle.
When understanding how the original variational structure is embedded at the discrete level, one can answer how the associated properties  will be preserved by the numerical solutions.
There has been a wide range of works about the conservation of geometrical properties 
at the numerical level by Hairer \textit{et al.}
\cite{HLW,MR2231943,MR3586395}, by Faou \cite{MR2895408} and on the conservation of variational structures by Marsden \textit{et al.} 
\cite{MR1462313,MR1697008,MW,MR2488601}
in the case of ODEs.
\\ \\
In this paper we will analyse the question of the conservation of variational structure as follows.
We consider the general framework of embeddings as presented in \cite{CD,cr1,CG,ci}.
We introduce the concept of \textit{coherence}.
Consider a problem associated to a Lagrangian structure. 
On one hand we have the Lagrangian functional $\LL$ on a functional space.
On the other hand we have the corresponding Euler-Lagrange equation.
Discretisation can be performed in two different ways.
\begin{itemize}
\item Either by discretising the Euler-Lagrange equation.
  This will be called a \textit{discrete differential embedding} because it
  is based 
  on deriving discrete versions of the differential operators in this PDE.
\item Or discretise the Lagrangian structure by defining a discrete Lagrangian functional $\LL_h$ and the associated discrete least action principle.
  This second procedure  is called \textit{discrete variational embedding} (it is also called {\sl variational integrator}).
\end{itemize}
In case the discrete differential embedding
and the discrete variational embedding are equivalent, we will say that we have \textit{coherence}.
% This is  enunciated in saying that the following diagram is commutative:
% \begin{eqnarray}
%   \nonumber
%   \begin{CD}
%     {{\rm Lagrangian }}~ \LL &  @>{{\rm disc.~ var.~ emb. }}>> &
%     {{\rm  discrete~ Lagrangian~ }} \LL_h\\
%     @V{{\rm L.A.P.}}VV & \hfill & @VV{{\rm  disc. ~L.A.P.}}V\\
%     {{\rm Euler-Lagrange~ equation}}  &  @>{{\rm disc.~ diff.~ emb.}}>> &   
%     {{\rm discrete~ Euler-Lagrange~ equation}}
%   \end{CD}
% \end{eqnarray}
% where
% L.A.P stands for  least-action principle. 
The same notion of coherence can be defined relatively to 
Hamiltonian structures.
\\
In case of coherence, the numerical solutions 
will inherit the properties of the original physical problem
(conservation of energy, independence with the coordinate system...).
% the discretisation firstly preserves the variational structure of the problem so transfering interesting properties (such as  independence with the coordinate system). It secondly may also preserve algebraic properties from the differential operators within the Euler-Lagrange equation.
\\
Based on this notion of coherence, the present work is an attempt to interpret 
numerical methods as {\it variational integrators} for PDEs deriving from  a Lagrangian/Hamiltonian structure. 
We will focus on a canonical example of such a problem: 
the Poisson equation.
This problem is well documented at the continuous  and at the discrete levels.
It provides an appropriate test case to improve the understanding of discrete embeddings for Lagrangian/Hamiltonian structure.
\\ \\
The outline of the paper is as follows. 
In section \ref{part1} are presented Lagrangian systems.
% In section \ref{sec:cl-euler-lag} the Lagrangian structure and the associated calculus of variations for fields are defined. The Lagrangian structure for the Poisson equation is recalled in section \ref{sec:lag-poisson}.
We introduce in section \ref{part:embedding}
% , the discrete embeddings are presented. A general definition of discrete embeddings is first presented in section \ref{generaldef}. D
the notions of discrete differential and discrete variational embeddings, 
and give various examples.
% are then defined in sections \ref{sec:disc-diff-emb} and \ref{sec:disc-var-emb} and illustrated with various examples. 
The concept of coherence % between differential and variational embeddings 
is then defined in section \ref{sec:coherence}.
% and illustrated by considering finite element methods.
In section \ref{part2}, we 
study the coherence for finite difference and finite volume methods, as applied to the Poisson equation.
% and on the coherence of two classical numerical methods for this problem:
% finite differences and finite volumes in sections \ref{section:fd} and \ref{section:fv} respectively.
Section \ref{part3} is concerned with Hamiltonian structures and mixed formulations.
% Hamiltonian structures and the associated calculus of variations are presented in section \ref{sec:ham}. 
% One recovers the mixed formulation of the Poisson equation with  Hamiltonian least action principle.
The discrete embedding of Hamiltonian structures 
is analysed for the mimetic finite difference method 
that is shown to be coherent.
% and the notion of 
% coherent embeddings are presented in section \ref{sec:disc-ham}.
%   % :   coherence  is shown to be naturally fulfilled by  conforming mixed finite element methods.
% Coherence of mimetic finite difference methods (see e.g. \cite{brezzi-2005})
%   % for the diffusion equation
% is analysed in the last section \ref{sec:mfd}.
\\ \\
Throughout this paper, $\Omega \subset \R^d $ is  a bounded domain with regular boundary.
The Sobolev space of order $m$ is denoted by ${\rm H}^m(\Omega)$ and the two following spaces
${\rm H}_0^1(\Omega)=\left\{ 
  v\in {\rm H}^1(\Omega), v_{|\partial\Omega}=0
\right\}$,     
${\rm H}_{\div}(\Omega)=\left\{ 
  \pp\in \left[{\rm L}^2(\Omega)\right]^d,~ 
  \div \pp \in {\rm L}^2(\Omega)
\right\}$
will be considered.

%%%%%%%%%%% New 10 juillet 2012%%%%%%%%%%%%%%%%%%
%%%%%%%%%%%%%%%%%%%%%%%%%%%%%%%%%%%%%%%%%%%%%%%% 
\section{Lagrangian systems}
%%%%%%%%%%%%%%%%%%%%%%%%%%%%%%%%%%%%%%%%%%%%%%%% 
%% 
%% 
\label{part1}
We recall classical results about Lagrangian calculus 
of variations for PDEs, illustrated in section \ref{sec:lag-poisson}  with the Lagrangian formulation of the Poisson problem.
For more details, we refer to \cite{evans,GH1,GH2}.
%% 
%% 
%% 
%% 
%% 
%%%%%%%%%%%%%%%%%%%%%%%%%%%% 
%%%%%%%%%%%%%%%%%%%%%%%%%%%% 
\subsection{Lagrangian calculus of variations}
%%%%%%%%%%%%%%%%%%%%%%%%%%%% 
%%%%%%%%%%%%%%%%%%%%%%%%%%%% 
%% 
%% 
\label{sec:cl-euler-lag}
\begin{definition}
  \label{defi:lagrange}
  An admissible Lagrangian function $L$ is a  function,
  \begin{eqnarray*}
    L :\Omega \times \R \times \R^d &\longrightarrow & \R\\
    (x,y,z)&  \mapsto &L(x,y,z),
  \end{eqnarray*}
  such that 
  $L$ is of class $\mathcal{C}^1$ with respect to $y$ and $z$ and integrable in 
  $x$.
  The Lagrangian function $L$ defines the \textit{Lagrangian functional} $\LL$:
  \begin{eqnarray*}
    \LL : {\rm H}^1(\Omega) &\rightarrow& \R, \\
    u & \longmapsto &
    \int_{\Omega} L(x,u(x),\nabla u(x))\,\dx.
  \end{eqnarray*}
\end{definition}
We are interested to vanish the first variations of the Lagrangian functional $\LL$ on a space of variations $V$. 
As in \cite{GH1}, we could give  a general notion for extremals
and variations. 
We take the following definitions of 
the notions of a differentiable functional and an extremal
for $\LL$.
\begin{definition}[Differentiability]
  \label{def:diff}
  We consider a space of variations $V\subset {\rm H}^1(\Omega)$. %$V\subset \DL$. 
  The functional $\LL$ is {\it differentiable} at point $u\in {\rm H}^1(\Omega) $ 
  % $u\in \DL $ 
  if and only if the limit,
  \begin{equation*}
    \nonumber
    \lim_{\epsilon \rightarrow 0} \frac{\LL(u+\epsilon v)-\LL(u)}{\epsilon},
  \end{equation*}
  exists in any direction $v \in V$.
  We then define the differential $D\LL(u)$ of $\LL$ at point $u$ as,
  \begin{equation*}
    v\in V\mapsto 
    D\LL (u)(v)=\lim_{\epsilon \rightarrow 0} \frac{\LL(u+\epsilon v)-
      \LL(u)}{\epsilon}.
  \end{equation*}
\end{definition}
With the  above definition of differentiability, one recovers 
the usual definition of the differential in case $V={\rm H}^1(\Omega)$ and $D\LL(u)$ 
is linear and continuous in $u$ on ${\rm H}^1(\Omega)$. The definition given here suffices to 
introduce extremals:
\begin{definition}[Extremals]
  \label{classicalextremal}
  A function  $u\in {\rm H}^1(\Omega)$ is an extremal for the functional 
  $\LL$ relatively to the space of variations $V\subset {\rm H}^1(\Omega)$ if
  $\LL$ is differentiable at point $u$ and:
  \begin{displaymath}
    \quad D\LL (u)(v)=0\,\quad {\rm for~any} \quad v\in V.
  \end{displaymath}
\end{definition}
\begin{proposition}
  \label{prop:diff-cond}
  If
  $x\mapsto \displaystyle{\frac{\partial L}{\partial y} (x,u(x),\nabla u(x))}$ 
  and $x\mapsto \displaystyle{\frac{\partial L}{\partial z}
    (x,u(x),\nabla u(x))}$ 
  respectively are in ${\rm L}^2(\Omega)$ and in $\left[{\rm L}^2(\Omega)\right]^d$, then the Lagrangian functional $\LL$ is differentiable at point $u\in {\rm H}^1(\Omega)$.
  \\
  In that case the differential is given for any $v\in {\rm H}^1(\Omega)$ by:
  \begin{equation}
    \label{eq:diff-L}
    D\LL(u)(v) = \int_\Omega 
    \left[
      \dfrac{\partial L }{\partial y} \bigl(x,u(x),\nabla u(x)\bigl)~v(x)
      +
      \dfrac{\partial L }{\partial z} \bigl(x,u(x),\nabla u(x)\bigl)\cdot\nabla v(x)
    \right]\dx.
  \end{equation}
\end{proposition}
\begin{proof}
  Using a Taylor expansion of $L$ at the point
  $(x,u+\epsilon v, \nabla(u+\epsilon v))$ in the variables $y$ and $z$
  leads to:
  \begin{equation*}
    L \bigl(x,u+\epsilon v, \nabla(u+\epsilon v)\bigl)  =  L (x,u, \nabla u)
    +\epsilon\,v\, \frac{\partial L}{\partial y}(x,u, \nabla u)
    +\nabla(\epsilon\,v)\cdot \frac{\partial L}{\partial z}(x,u, \nabla u)
    +o(\epsilon).
  \end{equation*}
  Integrating over the domain $\Omega$ gives:
  \begin{align*}
    \nonumber
    \LL(u+\epsilon v)  = \LL(u)
    &+\epsilon\int_{\Omega} v(x)\frac{\partial L}{\partial y}(x,u(x), \nabla u(x))\dx
    \\
    &+\epsilon\int_{\Omega}
    \nabla v(x)\cdot \frac{\partial L}{\partial z}(x,u(x), \nabla u(x)) \dx 
    +o(\epsilon),
  \end{align*}
  leading to \eqref{eq:diff-L}.
\end{proof}
Extremals of the functional $\LL$ can be characterised 
by an order $2$ PDE, called the {\it Euler-Lagrange equation}
given in the following theorem.
\begin{theorem}[Least action principle] \label{th:LAP}
  Consider a Lagrangian functional $\LL$  that
  satisfies the sufficient conditions of differentiability of
  proposition \ref{prop:diff-cond}
  at point $u\in {\rm H}^1(\Omega)$.
  Assume that $u$ is an extremal for a given space of variations $V$
  and that 
  $\frac{\partial L  }{ \partial z} (\cdot,u(\cdot),\nabla u(\cdot))\in {\rm H}_{\div}(\Omega)$.
  Moreover the subspace $V_0=\{v\in V, v=0 {\rm~on~} \partial\Omega\}$ is 
  supposed to be
  dense in ${\rm L}^2(\Omega)$. Then $u$ 
  satisfies the Euler-Lagrange equation:
  \begin{equation}
    \label{EL}
    \frac{\partial L}{\partial y}(x,u(x), \nabla u(x))
    -\div \left(
      \dfrac{\partial L}{\partial z}(x,u(x), \nabla u(x))
    \right)=0\,.
  \end{equation}
\end{theorem}
In the sequel we will denote $P$ the differential operator
associated to the Euler-Lagrange equation
  given by
  \begin{equation}
    \label{EL-P}
    P(u):=\frac{\partial L}{\partial y}(x,u(x), \nabla u(x))
    -\div \left(
      \dfrac{\partial L}{\partial z}(x,u(x), \nabla u(x))
    \right).
  \end{equation}
\begin{proof}
  Following \eqref{eq:diff-L} and using the Green 
  formula gives: $\forall ~v\in V_0$,
  \begin{align*}
    \int_\Omega 
    \left[
      \dfrac{\partial L }{\partial y} (x,u(x),\nabla u(x))
      -\div
      \left(
        \dfrac{\partial L }{\partial z} (x,u(x),\nabla u(x))
      \right)
    \right]v(x)
    \dx = 0,
  \end{align*}
  which implies \eqref{EL} by density of $V_0$ in ${\rm L}^2(\Omega)$. %${\rm L}^2(\Omega)$.
\end{proof}
%% 
%% 
%% 
%% 
%% 
%% 
%% 
%% 
%% 
%%%%%%%%%%%%%%%%%%%%%%%%%%%% 
%%%%%%%%%%%%%%%%%%%%%%%%%%%% 
\subsection{Lagrangian structure for the Poisson problem}
%%%%%%%%%%%%%%%%%%%%%%%%%%%% 
%%%%%%%%%%%%%%%%%%%%%%%%%%%% 
%% 
%% 
\label{sec:lag-poisson}
We consider the  Poisson problem on $\Omega$ 
for a homogeneous Dirichlet boundary condition:
find $u\in {\rm H}^2(\Omega)$,
\begin{equation}
  \label{eq:pb}
  -\Delta u =f \quad {\rm in}\quad   \Omega,
  \quad  {\rm and} \quad u=0 \quad {\rm on}
  \quad   \partial \Omega,
\end{equation}
for a data $f \in {\rm L}^2(\Omega)$. Assuming that $\Omega$ is bounded with a smooth boundary,  problem \eqref{eq:pb} has a unique solution $u\in {\rm H}^2(\Omega)\cap {\rm H}_0^1(\Omega)$.
\\
Equation \eqref{eq:pb} is
the differential formulation of  the Poisson problem.
Let us now present its  variational formulation.
We consider the Lagrangian function $ L $:
\begin{displaymath}
  L (x,y,z)=\frac{1}{2}  z\cdot z -f(x)y.
\end{displaymath}
The associated Lagrangian functional $\LL$ is given by,
\begin{equation}    
  \label{eq:L-Poisson}
  \LL(u)=\int_{\Omega} \big(~\dfrac{1}{2} \vert \nabla u \vert^2 ~-~ f u ~\big)~\dx\,. 
\end{equation}
%% and  is differentiable on ${\rm H}^1(\Omega)$.
%% 
%% 
%% 
%\begin{corollary}
  The differential formulation (\ref{eq:pb}) of the Poisson problem 
  is equivalent to, 
  % The Poisson problem under differential form (\ref{eq:pb}) has the following equivalent formulation, 
  \begin{equation}
    \label{eq:poisson-var}
    {\rm find }\quad  u\in {\rm H}^1_0(\Omega)\quad  {\rm so~that}\quad  \forall ~v\in {\rm H}^1_0(\Omega),\quad 
    D\LL(u)(v)=0,
  \end{equation}
  Equation \eqref{eq:poisson-var} is the well-known
  variational formulation of the 
  Poisson problem, with the space of variation $V={\rm H}^1_0(\Omega)$
  given by:
  \begin{equation*}
\int_{\Omega} \nabla u\,\nabla v\,\dx=\int_{\Omega} f\,v\dx\,.
    \end{equation*}
  .
%\end{corollary}
%% 
%% 
%% 
%% 
%% 
%% 
%% 
%% 
%%%%%%%%%%%%%%%%%%%%%%%%%%%%%%%% 
\section{Discrete embeddings}
\label{part:embedding}
%%%%%%%%%%%%%%%%%%%%%%%%%%%%%%% 5%

The formalism of embeddings has been initiated in \cite{CD} and 
further developed  in \cite{cr1,CG,ci}.
We propose here a general notion of  discrete embeddings.
% that applies to mappings on functional spaces.
This notion is defined in two particular cases: discrete embeddings of differential operators called discrete differential embedding in section \ref{sec:disc-diff-emb} and discrete  embeddings of Lagrangian functionals
called discrete variational embedding in
section \ref{sec:disc-var-emb}. 
The notion of coherence between discrete differential and discrete variational embeddings is presented in section \ref{sec:coherence}.

\subsection{General definitions} %NEW
\label{generaldef}
Let $X$ denote a
% some (infinite dimensional) 
functional space on $\Omega$. We consider the mapping,
$$
P:~u\in X \mapsto P(u) \in Y,
$$
where $Y$  either is a functional space on 
$\Omega$ or $Y=\R$. At this point no particular property is required for $P$.
\begin{definition}
  \label{def:disc-emb}
  We consider $X_h$ and $Y_h$ two finite dimensional spaces and 
  $\pi_1:~X\rightarrow X_h$, $\pi_2:~Y\rightarrow Y_h$ two 
  surjective linear mappings.
  We introduce $P_h:~X_h\rightarrow Y_h$ and consider the diagram:
  \begin{eqnarray}
    \label{diagramstar}
    \begin{CD}
      X &  @>{P}>> & Y\\
      @V{\pi_1}VV & \hfill & @VV{\pi_2}V\\
      {X_h}  &  @>{P_h}>> &   
      {Y_h}
    \end{CD}
  \end{eqnarray}
  We say that $P_h$ is a discrete embedding of $P$.
\end{definition}
\begin{remark}
  The setting presented in definition \ref{def:disc-emb} is  general.
  It introduces discrete (finite dimensional) counterparts for the functional spaces $X$ and $Y$.
  These discrete spaces theirselves can be functional spaces (such as for finite element methods \textit{e.g.}) or not (such as for finite difference methods).
  The diagram is not commutative in general.
\end{remark}
Consider again the Poisson problem. 
On one hand we have its differential formulation (\ref{eq:pb}). It is associated to the mapping $P:~u\in X \mapsto \Delta u +f\in Y$, with   $X={\rm H}^2(\Omega)$ and $Y={\rm L}^2(\Omega)$.
% Under its differential form (\ref{eq:pb}), it is associated with the mapping $P:~u\in X \rightarrow \Delta u +f\in Y$, with   $X={\rm H}^2(\Omega)$ and $Y={\rm L}^2(\Omega)$. 
The Poisson problem %, differential form, 
rewrites as:
\begin{displaymath}
  {\rm find}\quad  u\in M\subset X \quad {\rm so~ that}\quad 
  P(u) =0,
\end{displaymath}
with $M={\rm H}^1_0(\Omega)\cap X$.
A discretisation for the  differential formulation of the Poisson problem reads,
\begin{equation}
  \label{eq:poisson-Ph}
  {\rm find}\quad  u_h\in M_h\subset X_h \quad {\rm so~ that}\quad 
  P_h(u_h) =0,
\end{equation}
where $P_h:~X_h\rightarrow Y_h$ is a discrete embedding of $P$ and where $M_h\subset X_h$
% . In practice, the boundary condition either is directly included within the definition of $P_h$ in which case $M_h=X_h$, or is not included in the definition of $P_h$ in which case $M_h$ 
encodes the boundary condition.
The definition of $P_h$ requires a definition of $\Delta_h$.
This is a discrete embedding for the Laplace operator and will be referred as discrete differential embeddings.
This is detailed in section \ref{sec:disc-diff-emb}.
\\  \\
On the other hand the variational formulation  (\ref{eq:poisson-var})
of the Poisson problem, with $X={\rm H}^1(\Omega)$, 
$\LL:~X\rightarrow\R$ and $V={\rm H}^1_0(\Omega)=M$ rewrites as,
\begin{displaymath}
  {\rm find}\quad  u\in M\subset X \quad {\rm so~ that}\quad \forall ~v\in V,\quad 
  D\LL(u)(v)=0.
\end{displaymath}
A discretisation for the variational formulation of the Poisson problem reads:
\begin{displaymath}
  {\rm find}\quad  u_h\in M_h\subset X_h \quad {\rm so~ that}\quad \forall ~v_h\in V_h,\quad 
  D\LL_h(u_h)(v_h)=0.
\end{displaymath}
It involves $\LL_h:~X_h\rightarrow\R$, 
a discrete embedding of the Lagrangian functional $\LL:~X\rightarrow\R$,
that will be referred as discrete variational embedding.
This is developed in section \ref{sec:disc-var-emb}.
\subsection{Discrete differential embeddings}
\label{sec:disc-diff-emb}
\begin{definition}%[Discrete differential embedding]
  \label{def:disc-diff-emb}
  Consider the diagram (\ref{diagramstar}) in definition \ref{def:disc-emb}
  in the case where $P$ is associated with some PDE $P(u)=0$, \textit{i.e.}  $P$ is a differential operator.
  In that particular case we call 
  $P_h$  a discrete differential embedding.
\end{definition}
Note that a discrete differential embedding is not a differential operator itself. It is the discretisation of a differential operator.
\\ \\
% \begin{remark}
Consider the discrete differential embedding for the Poisson problem
\eqref{eq:poisson-Ph}.
We set $P_h u_h= \Delta_h u_h+ f_h$.
The definition of $P_h$ involves a definition of $f_h$ and of $\Delta_h$. 
Two ways can be followed to derive $\Delta_h$. The first one is to directly discretise the Laplacian, as it is done using finite difference methods in section \ref{section:fd}. The second one is to use the divergence form of the Laplacian: $\Delta= \div\circ \nabla$ and to derive a discrete embedding for the Laplacian as $\Delta_h= \div_h\circ \nabla_h$,
where $\div_h$ and $\nabla_h$ are two discrete
differential embedding of $\div$ and $\nabla$. This will be the case with finite volume methods in section \ref{section:fv}.
\\
This leads to two discrete differential embeddings for the Poisson problem: either,
\begin{displaymath}
  -\Delta_h u_h =f_h,
\end{displaymath}
or,
\begin{displaymath}
  -\div_h (\nabla_h u_h) = f_h .
\end{displaymath}
These two discrete problems do not coincide in general. 
Indeed, recovering the  algebraic properties of the original 
differential operators (here $\Delta=\div \circ \nabla$)  
at the discrete level (here $\Delta_h=\div_h \circ \nabla_h$) 
is a full problem by itself.   
% \end{remark}
\\ \\
We now give three illustrations of discrete differential embeddings: 
for the gradient operator and for the divergence one.
Let us start precising the notion of a mesh for the domain 
$\Omega\subset \R^d$ $d=2,3$. 
\begin{definition}[Mesh]%[Mesh $\msh$, cells $ K $, faces $\fc$ and vertexes $\vtc$]
  \label{def:msh}
  A cell is a polygonal/polyhedral non empty open subset
  of $\R^d$.
  A mesh $\msh$ of the domain $\Omega$ is a collection of cells 
  partitioning $\Omega$ in the following sense:
  \begin{align*}
    \cup_{ K \in\msh} \overline{ K }=\overline{\Omega},\qquad \text{ and }
    \quad
    \Bigl(  K_1, K_2\in\msh~\Rightarrow
    \quad {\rm either}\quad  K_1\cap K_2=\emptyset \quad {\rm or}\quad  K_1= K_2\Bigl).
  \end{align*}
  A face (or an edge) $ e $ of some $ K \in\msh$ such that $ e \subset \partial\Omega$ is called a boundary face. 
  The set of boundary faces is denoted  $\bfc$. It satisfies:
  $\partial\Omega=\displaystyle{\cup_{ e \in\bfc} e}$. 
  For every $ e \in\bfc$, there exists a unique $ K \in\msh$ satisfying
  $ e \subset \overline{ K }\cap\partial\Omega$: one writes $ e = K \vert \partial\Omega$. 
  \\
  The internal faces set $\ifc$ associated with $\msh$ is 
  the set of all 
  geometrical subsets $ e =\overline{ K_1}\cap\overline{ K_2}$,
  $ K_1, K_2\in\msh$ and $ K_1\neq K_2$, having non-zero
  $(d-1)-$dimensional measure. 
  For every $ e \in\ifc$, there exists a unique couple
  $ K_1, K_2\in\msh$ satisfying 
  $ e =\overline{ K_1}\cap\overline{ K_2}$: one  writes
  $ e = K_1\vert  K_2$. 
  \\
  The faces set associated with $\msh$ is given as 
  $\fc=\bfc\cup\ifc$. It provides a  partitioning of 
  $\displaystyle{\cup_{ K \in\msh}
    \partial  K }$, in the same meaning as earlier: $\displaystyle{\cup_{ e \in\fc}
    e }=\displaystyle{\cup_{ K \in\msh} \partial  K }$ 
  and the overlapping of two distinct
  faces either is  empty or of zero $(d-1)-$dimensional measure.
  Let $ e \in\fc$ such that $ e \subset \partial K $ for $ K \in\msh$.
  We denote $\bn{_{K, e }}$ the unit normal to $ e $ pointing outward of $ K $.
  We also provide intrinsic orientation to faces: to all faces $e\in\fc$ is associated $\bn_e$ one of its (two) unit normal, if $ e \subset \partial K $ we have $\bn_e=\pm \bn{_{K, e }}$.
  \\
  The set of vertexes associated with $\msh$ is denoted
  $\vtc$: it contains exactly all the vertexes of all the cells
  $ K \in\msh$. 
  \\ \\
  One shall denote $\meas{O}$ the measure of a geometrical object
  $O$ according to its dimension. Taking $d=3$, $\meas{ K }$ is the
  volume of the cell $ K $, $\meas{ e }$ the area of an edge $ e \in\fc$ and
  $\meas{xy}$ the length between two points $x$ and $y$.
  The cardinal of a set $E$ is $\# E$. 
\end{definition}

\subsubsection{The finite volume divergence}
\label{sec:fv-div}
We denote here 
$X=[{\rm H}^1(\Omega)]^d$, $Y={\rm L}^2(\Omega)$ 
% $X={\rm H}_{\div}(\Omega)$, $Y={\rm L}^2(\Omega)$ 
and $\div:~X\rightarrow  Y$ is the divergence operator. 
Let $\msh$ be a mesh of $\Omega$. We here define $X_h = \R^{\# \fc}$,
% with $\# \fc$ the total number of faces of the mesh
and $Y_h=P^0(\msh)$ the space of piecewise constant functions over the cells of the mesh, with the natural identification $Y_h= \R^{\# \msh}$.
% with $\# \msh$ the total number of cells of the mesh.
Note that in general there is no natural identification of $ \R^{\# \fc}$ with some finite dimensional vector field space over $\Omega$, we however mention the case of simplicial meshes where such an identification is provided by the Raviart-Thomas finite element space of order 0, $RT_0(\Omega)$, see \cite{RT77}. 
\\ 
To $\pp\in X$ we associate $\pi_1 \pp=(p_e)_{e\in \fc}$ with $p_e=\int_e \pp\cdot \bn_e \dl / |e|$
the mean flux of $\pp$ across the face $e$ according to its orientation provided by $\bn_e$ (in the trace sense). 
To $f\in {\rm L}^2(\Omega)$, we associate $\pi_2 f=(f_K)_{K\in \msh}$ with $f_K=\int_K f\, \dx / |K|$ the mean value of $f$ on the cell $K$.
The discrete divergence is defined as,
$$
\div_h:~ \pp_h =(p_e)_{e\in\fc} \in  \R^{\# \fc} \mapsto (\div_K \pp_h)_{K\in \msh} \in \R^{\#\msh},
$$
with,
\begin{equation}
  \label{eq:divh-fv}
  \div_K \pp_h = \dfrac{1}{|K|} \sum_{e\in\fc, e \subset \partial K} p_e |e| ~\bn_e\cdot \bn{_{K, e }}.    
\end{equation}
This definition simply is the flux balance around the cell $K$, the last term $\bn_e\cdot \bn{_{K, e }}$ giving the correct orientation for the fluxes, 
{\it i.e.} outside the cell $K$.
\\ \\
With these definitions we have a discrete differential embedding for the divergence,
\begin{eqnarray}
  \label{eq:disc-div-fv}
  \begin{CD}
    [{\rm H}^1(\Omega)]^d &  @>\div>> & {\rm L}^2(\Omega)\\
    @V{\pi_1}VV & \hfill & @VV{\pi_2}V\\
    \R^{\# \fc}  &  @>\div_h>> &   
    R^{\#\msh}
  \end{CD}
\end{eqnarray}
and this diagram moreover is commutative thanks to the 
divergence formula: $\pi_2\circ\div = \div_h\circ\pi_1$.
\subsubsection{The $P^1(\msh)$ finite element gradient }
\label{sec:fe-grad}
We introduce $X=C^1(\Omega)$ and $Y=\left[C^0(\Omega))\right]^d$  the spaces of continuously differentiable functions and of continuous vector fields  over $\Omega$ respectively. We now consider the gradient operator $\nabla :~ C^1(\Omega) \rightarrow \left[C^0(\Omega)\right]^d$. 
\\
Let $X_h = P^1(\msh)$ be the space of continuous functions over $\Omega$ that moreover are piecewise affine on each cell $K\in\msh$.
Let us assume that the mesh is simplicial, the space $X_h$ is identified 
to $\R^{\# \vtc}$.
% with $\#\vtc$ the total number of vertexes of the mesh $\msh$.
We have the projection $\pi_1:~u\in C^1(\Omega)\mapsto \pi_1 u=(u_S)_{S\in\vtc}\in  P^1(\msh)$ with $u_S=u(S)$.
Let $Y_h=\left[ P^0(\msh) \right]^d$ be the space of piecewise constant vector fields over each cell $K\in\msh$. We  have a simple projection $\pi_2:~\left[C^0(\Omega)\right]^d\rightarrow \left[ P^0(\msh) \right]^d$ by averaging a vector field over each cell of the mesh (similarly to $\pi_2$ in section \ref{sec:fv-div}).
\\ \\
We have the following discrete differential embedding for the gradient:
\begin{eqnarray}
  \label{eq:disc-grad-P1}
  \begin{CD}
    C^1(\Omega)&  @>\nabla>> & \left[C^0(\Omega))\right]^d\\
    @V{\pi_1}VV & \hfill & @VV{\pi_2}V\\
    P^1(\msh)  &  @>\nabla _h>> &   
    \left[ P^0(\msh) \right]^d
  \end{CD}
\end{eqnarray}
where the discrete gradient $\nabla _h = \nabla _{| P^1(\msh)}$ indeed is the restriction of 
the continuous one to $P^1(\msh)$.
In that case the diagram is not commutative.

\subsubsection{Non-conforming finite element gradient }
We  propose a second definition of discrete differential embedding of 
the gradient, that is referred to as \textit{non-conforming finite 
  element gradient} since it matches with the Crouzeix-Raviart finite element 
of order 1 discretisation,  see
\cite{CR73}, in the case of a simplicial mesh.
\\
Let $X={\rm H}^1(\Omega)$, $Y=\left[{\rm L}^2(\Omega)\right]^d$ and consider 
$\nabla :~{\rm H}^1(\Omega) \rightarrow \left[{\rm L}^2(\Omega)\right]^d$.
We set $X_h = \R^{\# \fc}$
% , with $\# \fc$ the total number of faces of the mesh
and $Y_h=\left[P^0(\msh)\right]^d$ the space of piecewise constant vector fields over the cells of the mesh, with the natural identification $Y_h= \left[\R^d\right]^{\# \msh}$.
% , with $\# \msh$ the total number of cells of the mesh.
% 
To $u\in {\rm H}^1(\Omega)$ we associate $\pi_1 u=(u_e)_{e\in \fc}$ with $u_e=\int_e u \dl / |e|$ the mean value of $u$ on the face $e$ (in the trace sense).
We  have the same simple projection $\pi_2:~\left[{\rm L}^2(\Omega)\right]^d\rightarrow \left[\R^d\right]^{\# \msh}$ as in
section \ref{sec:fe-grad}  by averaging a vector field over each cell of the mesh.
The discrete gradient is defined as,
$$
\nabla _h:~ u_h =(u_e)_{e\in\fc} \in X_h \mapsto (\nabla _K u_h)_{K\in \msh} \in \left[\R^d\right]^{\# \msh},
$$
with,
$$
\nabla _K u_h = \dfrac{1}{|K|} \sum_{e\in\fc, e \subset \partial K} u_e |e| ~ \bn{_{K, e }}.
$$
With these definitions we have the following discrete differential embedding for the gradient,
\begin{eqnarray*}
  \label{eq:disc-grad-nc}
  \begin{CD}
    {\rm H}^1(\Omega) &  @>\nabla >> & \left[{\rm L}^2(\Omega)\right]^d\\
    @V{\pi_1}VV & \hfill & @VV{\pi_2}V\\
    X_h  &  @>\nabla _h>> &   
    \left[\R^d\right]^{\# \msh}
  \end{CD}
\end{eqnarray*}
and this diagram moreover is commutative thanks to the  formula $\int_K \nabla u \dx = \int_{\partial K} u \bn \dl$, with $\bn$ the unit normal on $\partial K$ pointing outwards $K$.
\subsection{Discrete variational embeddings}
\label{sec:disc-var-emb}
\begin{definition}%[Discrete variational embedding]
  \label{def:disc-var-emb}
  We consider a Lagrangian functional $\LL:~X\rightarrow \R$ as defined in definition \ref{defi:lagrange} for some 
  % (infinite dimensional) 
  functional space $X\subset {\rm H}^1(\Omega)$. 
  A discrete variational embedding is a discrete embedding $\LL_h$ of $\LL$ as defined in definition \ref{def:disc-emb} in the particular framework $Y=\R=Y_h$ and $\pi_2=id$. The diagram for a discrete variational embedding is the following,
  % \centerline{
  $$
  \xymatrix{
    X  \ar[r]^{\LL} \ar[d]_{\pi_1}  & \R  \\
    X_h \ar[ur]_{\LL_h} & 
  }
  $$
  % }
\end{definition}

\subsubsection*{ Finite element discrete variational embedding}
% \label{sec:fe-var-emb}
We consider a general Lagrangian functional $\LL$ as
in definition \ref{defi:lagrange}.
We use here the same framework as in section \ref{sec:fe-grad}. The mesh is assumed to be simplicial. We consider $X=C^1(\Omega)$, $X_h=P^1(\msh)$ 
and the projection 
$\pi_1:~u\in C^1(\Omega)\mapsto \pi_1 u=(u_S)_{S\in\vtc}\in X_h$ 
with $u_S=u(S)$. 
Since $P^1(\msh)\subset {\rm H}^1(\Omega)$ we define $\LL_h:~P^1(\msh)\rightarrow\R$ as $\LL_h=\LL_{|P^1(\msh)}$. 
We have the diagram,
% \centerline{
$$
\xymatrix{
  C^1(\Omega)  \ar[r]^{\LL} \ar[d]_{\pi_1}  & \R  \\
  P^1(\msh) \ar[ur]_{\LL_h} & 
}
$$
% }
Note that this definition  extends to any conformal finite
element space $X_h$, since we always have $X_h\subset {\rm H}^1(\Omega)$ (see e.g.\cite{Ciarlet,GR}).
Of course the definition of $\pi_1$ needs to be adapted to each particular choice of $X_h$.
\\
Also note that the extension to non-conforming finite elements is possible since $\LL$ can be evaluated on any function $u$ that would only be locally ${\rm H}^1$, over each cell of the mesh (precisely $u_{|K}\in {\rm H}^1(K)$ for all $K\in\msh$) instead than globally ${\rm H}^1$ on the whole domain $\Omega$.
%% 
%% 
%% 
%% 
%%%%%%%%%%%%%%%%%%%%%%%%%%%% 
%%%%%%%%%%%%%%%%%%%%%%%%%%%% 
\subsection{Coherence}
%%%%%%%%%%%%%%%%%%%%%%%%%%%% 
%%%%%%%%%%%%%%%%%%%%%%%%%%%% 
\label{sec:coherence}
Consider a problem associated with a Lagrangian variational structure and consider this problem either under its variational formulation 
(Lagrangian least action principle),
\begin{equation}
  \label{eq:pde-var}
  {\rm find}\quad  u\in M\subset X \quad {\rm so~ that}\quad \forall ~v\in V,\quad 
  D\LL(u)(v)=0,
\end{equation}
or under its differential formulation (Euler-Lagrange equation),
\begin{equation}
  \label{eq:pde-diff}
  {\rm find}\quad  u\in M'\subset X' \quad {\rm so~ that}\quad 
  P(u) =0,
\end{equation}
where $P(u)$ defined in equation \eqref{EL-P}
is the operator associated to the Euler-Lagrange equation.
\\
Under the conditions of theorem \ref{th:LAP}, these two formulations are equivalent. 
They however give rise to two discretisation procedures.
\begin{itemize}
\item [--] Being given $\LL_h$ a discrete variational embedding of $\LL$ as in definition \ref{def:disc-var-emb}, the discrete least action principle reads
  \begin{equation}
    \label{eq:pdeh-var}
    {\rm find}\quad  u_h\in M_h\subset X_h \quad {\rm so~ that}
    \quad \forall ~v_h\in V_h, \quad
    D\LL_h(u_h)(v_h)=0.
  \end{equation}
  This is a discrete variational formulation.
\item [--] Being given a discrete differential embedding of $P$ as in definition \ref{def:disc-diff-emb}, 
  \begin{equation}
    \label{eq:pdeh-diff}
    {\rm find}\quad  u_h\in M_h'\subset X_h' \quad {\rm so~ that}\quad 
    P_h(u_h) =0.
  \end{equation}
  This is a discrete differential formulation. 
\end{itemize}
\textit{A priori,} the two discrete problems \eqref{eq:pdeh-var} and \eqref{eq:pdeh-diff} do not  provide equivalent problems.
This question is addressed considering the concept of \textit{coherence} 
introduced in \cite{CD}. 
\begin{definition}[Coherence]
  \label{def:coherence}
  Consider a Lagrangian functional $\LL$ satisfying the hypothesis of theorem 
  \ref{th:LAP}.
  The operator associated to its Euler Lagrange equation  \eqref{EL} is denoted $P$.
  \\
    Consider a discrete variational embedding of $\LL$ as in definition \ref{def:disc-var-emb}.
    It is associated to a functional $\LL_h$
    and to a discrete least action principle given by equation \eqref{eq:pdeh-var}.
  Consider 
  a discrete differential embedding of $P$ as in definition \ref{def:disc-diff-emb}:
  it is associated to an operator $P_h$.
  These embeddings are said coherent if the discrete variational formulation
  \eqref{eq:pdeh-var}
  and the discrete differential formulation 
  \eqref{eq:pdeh-diff}  are equivalent (\textit{i.e.} have the same solutions).
  \\
  In other words the following  diagram is commutative, 
  \begin{eqnarray*}
    \begin{CD}
      {{u \mapsto \LL(u)}} &   @>{{\rm disc.~var.~emb.}}>> & 
      { {u_h \mapsto \LL_h(u_h)}}\\
      @V{{\rm L.A.P.}}VV & \hfill & @VV{ {\rm disc.~L.A.P.}}V\\
      u \;{\rm solution~of~PDE}~\eqref{eq:pde-diff}  & @>{{\rm disc.~diff.~emb.}}>> &   
      u_h {\rm ~solution~of~PDE}_h~ \eqref{eq:pdeh-diff}\\
      {\rm E.L.~equation  }  & @.& 
      {{\rm disc.~E.L.~equation }}
    \end{CD}
  \end{eqnarray*}
  where L.A.P. stands for least action principle and E.L. for Euler Lagrange.
\end{definition}
A general raised  question then is:
can we find conditions ensuring  the coherence between the discrete differential and variational embeddings ? 
\\
In the next two parts we study the coherence for discrete differential embeddings of problems having a Lagrangian or Hamiltonian variational formulation.
It turns out that one cannot set apart the coherence from the algebraic properties of $P_h$ inherited from the one of $P$. More precisely a property of integration by parts type is required at the discrete level to ensure coherence.
\\
A deeper insight into this relationship is gained by considering the Poisson problem. Assume one  performs a discrete differential embedding $\Delta_h$ for the Laplacian.
In all forthcoming examples,
coherence is obtained in case  $\Delta_h$  is the composition of a discrete gradient and a discrete divergence $\Delta_h=\div_h\circ\nabla _h$, and if in addition these two discrete operators fulfil a duality property of type 
Green-Gauss formula.
This is the case for finite differences with formula \eqref{eq:DF-green}, for  finite volumes with formula \eqref{eq:fv-green} and for mimetic finite differences in section \ref{sec:mfd}.
\subsubsection*{Coherence for conforming finite elements}
We first consider the 
% Euler-Lagrange PDE \eqref{EL} 
Poisson problem (\ref{eq:pb}).
As developed in section \ref{sec:lag-poisson}, this PDE is the Euler Lagrange equation associated with a least action principle on the Lagrangian functional $\LL(u)=\int_{\Omega} \big(\dfrac{1}{2} \vert \nabla u \vert^2 - f u \big)\dx$  
given in equation \eqref{eq:L-Poisson}.
\\
Let $X_h\subset  {\rm H}_0^1(\Omega)$ be some conforming finite element space.
We can define $\LL_h=\LL_{\vert  X_h}$. 
This provides a discrete variational embedding of $\LL$ as in definition \ref{def:disc-var-emb}.
The numerical problem solved in practice is the linear problem 
$P_h(u_h)=0$   on $X_h$ 
(involving the mass and stiffness matrices) where
the operator $P_h$ is defined by,
\begin{displaymath}
  \forall v_h\in X_h,\quad
  \int_\Omega P_h(u_h)v_h \dx = 
  \int_\Omega \left (\nabla u_h \cdot  \nabla v_h - f v_h\right )\, \dx
  =D\LL_h(u_h)(v_h).
\end{displaymath}
The operator $P_h$ on $X_h$ provides a discrete differential embedding for the operator $P(u) = -\Delta u -f$ but is not explicit.
By construction, these two discrete variational and differential embeddings are coherent.
\\ \\
The coherence for conforming finite element methods naturally extends to the PDE $P(u)=0$ in equation (\ref{EL}) for an homogeneous Dirichlet boundary condition. 
This problem derives from a least action principle associated with the Lagrangian functional $\LL$ in definition \ref{defi:lagrange}. On one hand we have a discrete variational embedding with $\LL_h=\LL_{\vert  X_h}$.
On the other hand the problem solved in practice is $P_h(u_h)=0$ with $P_h(u_h)$ defined as,
\begin{equation*}
  % \label{eq:fe-disc}
  \forall ~ v_h\in X_h,\quad
  \int_\Omega P_h(u_h)\,v_h\,\dx=
  \int_\Omega 
  \Bigl(
  \dfrac{\partial L }{\partial y} (x,u_h,\nabla u_h)\,v_h
  +
  \dfrac{\partial L }{\partial z} (x,u_h,\nabla u_h)\cdot\nabla v_h
  \Bigl)\,\dx,
\end{equation*}
that provides a discrete differential embedding of $P$.
These discrete embeddings are coherent by construction.
%% 
%% 
%% 
%% 
%% 
%% 
%% 
%% 
%% 
%% 
%%%%%%%%%%%%%%%%%%%%%%%%%%%%%%%%%%%%%%%%%%%%%%%% 
%%%%%%%%%%%%%%%%%%%%%%%%%%%%%%%%%%%%%%%%%%%%%%%% 
\section{Coherence of classical discrete embeddings}
%%%%%%%%%%%%%%%%%%%%%%%%%%%%%%%%%%%%%%%%%%%%%%%% 
%%%%%%%%%%%%%%%%%%%%%%%%%%%%%%%%%%%%%%%%%%%%%%%% 
%% 
%% 
\label{part2}
In section \ref{sec:coherence},  we showed a first example of coherent discrete embedding of Lagrangian structure.
In this precise case, several facilities were available: the discrete solution also is  a function $u_h:~\Omega\rightarrow \R$ so that differentiation and integration had the same sense at the discrete and at the continuous levels.
As a result the definition of a discrete Lagrangian $\LL_h$ was obvious 
and natural: $\LL_h$ 
was the restriction of $\LL$ to some functional space of finite dimension.
\\
Such facilities are not always available, they rather are restricted to conforming finite element methods. Such a lifting between the discrete space of unknowns $X_h$  and a function space is not available in general.
As a result differentiation and integration have to be re-defined at the discrete level to provide a definition of a discrete Lagrangian. In this section we give two examples of discrete embeddings for a Lagrangian structure: finite differences and classical finite volumes. 
Coherence is proved in both cases.
%% 
%% 
%% 
%% 
%% 
%%%%%%%%%%%%%%%%%%%%%%%%%%%% 
%%%%%%%%%%%%%%%%%%%%%%%%%%%% 
\subsection{Finite differences}
%%%%%%%%%%%%%%%%%%%%%%%%%%%% 
%%%%%%%%%%%%%%%%%%%%%%%%%%%% 
%% 
%% 
%% 
%% 
%% 
\label{section:fd}
We refer to \cite{thomee} for a general presentation of finite difference methods.
We study in this section the coherence properties of finite difference methods applied firstly to the Poisson problem (\ref{eq:pb})  and secondly to the general Euler-Lagrange PDE \eqref{EL}.
The domain is set to $\Omega=[0,1]^2$. We consider a 
Cartesian grid $\msh$ of $\Omega$ with uniform size $h=1/N$, $N\in\N^*$,
in every direction. The results of this section can be extended 
to more general domains, to other space dimensions and more general lattices.
\\
We will use the following notations.
For 
$\bj=(i,j) \in \Z^2$ we write $0 \leq \bj\leq N$ 
if $ 0\leq i\leq N$ and $ 0\leq j\leq N$. Let
$J=\{\bj\in\N^2,~0\le\bj\le N\}$. 
The point of coordinates $(i h,j h)\in \overline{\Omega}$ is denoted $x_{\bj}$.
The mesh  with vertexes 
$\{x_\bj,~\bj\in J \}$ is denoted $\msh$, it is a cartesian grid of $\Omega$.
\\
Let us consider the two spaces
$\mathcal{S}=\{u:~\Z^2\longrightarrow \R\}$ and 
$\mathcal{V}=\{\pp:~\Z^2\longrightarrow \R^2\}$.
Let $\bj=(i,j)\in\Z^2$: we denote for $u\in\mathcal{S}$,  $u_{i,j} = u_\bj = u(\bj)$ and for $\pp\in\mathcal{V}$,
$\pp_{i,j} = \pp_\bj = \pp(\bj)$.
We consider the discrete  space
$$
X_h=\{u\in\mathcal{S},~u_\bj=0~~{\rm if}~~\bj\notin J ~~{\rm and~ if}~~x_\bj\in \partial \Omega\}.
$$
The truncation operator $T:~\mathcal{S}\rightarrow X_h$ is defined as
$\bigl(Tu\bigl)_\bj=u_\bj$ if $0<\bj<N$ and by $\bigl(Tu\bigl)_\bj=0$ otherwise.  
\subsubsection{Discrete differential embedding for the Laplacian}
\label{sec:df-diff-emb}
The discrete Laplacian $\Delta_h:~\mathcal{S} \rightarrow \mathcal{S}$ is defined as, for $\bj=(i,j)\in\Z^2$,
$$
(\Delta_h u)_\bj = 
\dfrac{u_{i-1,j}-2 u_{i,j} + u_{i+1,j} }{h^2}
+
\dfrac{u_{i,j-1}-2 u_{i,j} + u_{i,j+1} }{h^2}.
$$
The operator  $T\circ\Delta_{h}:~\mathcal{S}\rightarrow X_h$ induces a mapping on $X_h$.
Considering the projection $\pi_1:~C^0(\Omega)\rightarrow X_h$,
given by $(\pi_1 u)_\bj = u(x_\bj)$ if $0<\bj<N$, or $(\pi_1 u)_\bj=0$ otherwise, we have the following discrete differential embedding:
\begin{eqnarray*}
  \begin{CD}
    C^2(\Omega) &  @>\Delta>> & C^0(\Omega)\\
    @V{\pi_1}VV & \hfill & @VV{\pi_1}V\\
    X_h  &  @>T\circ\,\Delta_h>> &   
    X_h
  \end{CD}
\end{eqnarray*}
For $f\in C^0(\Omega)$, the discrete differential embedding of $P(u)=-\Delta u -f$
then is $P_h(u_h) = -T\circ\Delta_h u_h -\pi_1 f$ for $u_h\in X_h$.
The discrete differential formulation of 
the Poisson problem is,
\begin{equation}
  \label{eq:df-df}
  {\rm find }\quad u\in X_h\quad {\rm so~ that} \quad 
  P_h(u) = -T\circ\,\Delta_h u - \pi_1 f = 0.  
\end{equation}
\\ 
Let us introduce a discrete gradient 
$\nabla_h:~\mathcal{S} \rightarrow \mathcal{V}$
and a discrete divergence 
$\div_h:~ \mathcal{V}\rightarrow \mathcal{S}$.
For $\bj=(i,j)\in\Z^2$ they are given by,
\begin{displaymath}
  % \label{eq:def-grad-div-df}
  (\nabla_h u)_\bj= \dfrac{1}{h}\left(
    \begin{array}{c}
      u_{i+1,j}-u_{i,j}       
      \\[5pt]
      u_{i,j+1}-u_{i,j}
    \end{array}
  \right),\quad 
  (\div_h \pp)_\bj= 
  \dfrac{p^1_{i,j}     - p^1_{i-1,j}}{h} + 
  \dfrac{p^2_{i,j}     - p^2_{i,j-1}}{h}
  ,
\end{displaymath}
for $u\in \mathcal{S}$ and $\pp=(p^1, p^2)\in\mathcal{V}$ ($p^1\in\mathcal{S}$ and $p^2\in \mathcal{S}$ are the two components of $\pp$).
This defines two discrete differential embeddings,
\begin{eqnarray*}
  \begin{CD}
    C^1(\Omega) &  @>\nabla >> & \left[C^0(\Omega)\right]^2\\
    @V{\pi_1}VV & \hfill & @VV{\pi_2}V\\
    X_h  &  @>T_2\circ\,\nabla _h>> &   
    X_h\times  X_h
  \end{CD}\quad \quad,\qquad \qquad 
  \begin{CD}
    \left[C^1(\Omega)\right]^2 &  @>\div >> & C^0(\Omega)\\
    @V{\pi_2}VV & \hfill & @VV{\pi_1}V\\
    X_h\times  X_h  &  @>T\circ\,\div _h>> &   
    X_h
  \end{CD}
\end{eqnarray*}
with $\pi_2=\pi_1\times \pi_1$ and  $T_2=T\times T$ component by component.
\\
As one can see, a forward finite difference formula has been used for the definition of the discrete gradient, whereas a backward one has been used for the discrete divergence.
This choice has been made in order to have the following properties
(that can easily be checked).
We have the composition rule,
\begin{equation}
  \label{eq:chain-rule-df}
  \Delta_h = \div_h \circ \nabla_h,
\end{equation}
and the \textit{discrete Green-Gauss formula}:
\begin{equation}
  \label{eq:DF-green}
  \forall ~u\in X_h,\quad  \forall ~\pp\in\mathcal{V}:\quad 
  \sum_{ \bj\in J}
  \pp_\bj\cdot(\nabla_h u)_\bj
  =-
  \sum_{ \bj\in J}
  (\div_h \pp)_\bj ~ u_\bj.
\end{equation}
\subsubsection{Discrete variational embedding, coherence}
For $f\in C^0(\Omega)$  we introduce the discrete Lagrangian functional:
\begin{displaymath}
  \forall ~u\in X_h,\quad \LL_h(u) = 
  \dfrac{1}{2}\sum_{\bj\in J} |\nabla _{\bj} u|^2 ~h^2 ~-~
  \sum_{\bj\in J} (\pi_1 f)_{\bj} u_{\bj}~h^2.
\end{displaymath}
This definition provides the following discrete variational embedding for the Poisson Lagrangian functional 
$\LL:~u\mapsto \int_\Omega(\frac{1}{2}|\nabla u|^2-fu)\dx$ given
in equation (\ref{eq:L-Poisson}),
% This definition provides the following discrete variational embedding gor the Poisson Lagrangian functional $\LL:~u\longrightarrow \int_\Omega(|\nabla u|^2/2-fu)\dx$,
$$
\xymatrix{
  C^1(\Omega)  \ar[r]^{ \LL} \ar[d]_{\pi_1}  & \R  \\
  X_h \ar[ur]_{\LL_h} & 
}
$$
The discrete variational formulation of the Poisson problem reads:
\begin{equation}
  \label{eq:DF-ELh}
  {\rm find }\, u\in X_h \quad {\rm so~ that}\quad 
  \forall ~v_h\in X_h:\quad 
  D\LL_h(u)(v_h)=0.
\end{equation}
\begin{theorem}
  \label{thm:coherence-df}
  The discrete variational and differential embeddings of the Poisson problem using  the finite difference method  are coherent. 
  Precisely, the two discrete problems \eqref{eq:df-df} and \eqref{eq:DF-ELh} have the same solutions.
\end{theorem}
\begin{proof}
  Let us consider a solution $u$ to \eqref{eq:DF-ELh}. We  have for all $v_h\in X_h$,
  \begin{displaymath}
    \sum_{\bj\in J} (\nabla_h u)_\bj\cdot(\nabla_h v_h)_\bj ~h^2
    -
    \sum_{\bj\in J} (\pi_1 f)_\bj v_\bj~h^2 =0.
  \end{displaymath}
  Using the discrete Green-Gauss formula \eqref{eq:DF-green}, we get: for all $v_h\in X_h$,
  \begin{displaymath}
    - \sum_{\bj\in J} \left (\div_h (\nabla_h u)\right )_\bj ~ v_\bj ~h^2
    -
    \sum_{\bj\in J} (\pi_1 f)_\bj v_\bj~h^2 =0.
  \end{displaymath}
  Using the composition rule \eqref{eq:chain-rule-df}, this exactly means, for all $\bj$ so that $0<\bj<N$,
  $
  -\left (\Delta_h u\right )_\bj = (\pi_1 f)_\bj,
  $
  which is equation \eqref{eq:df-df}.
\end{proof}
\subsubsection{Extension}
The previous coherence theorem extends to the general Euler-Lagrange PDE \eqref{EL} that we recall,
\begin{displaymath}
  P(u)=\frac{\partial L}{\partial y}(x,u(x), \nabla u(x))
  -\div \left(
    \dfrac{\partial L}{\partial z}(x,u(x), \nabla u(x))
  \right)=0.
\end{displaymath}
This equation is considered here together with a homogeneous boundary condition on $\partial\Omega$.
The two discrete differential embeddings for the gradient and for the divergence introduced in section \ref{sec:df-diff-emb} provide the following 
discrete differential embedding $P_h:~X_h\rightarrow X_h$. It is defined for $u\in X_h$ by,
\begin{align*}
  % \label{eq:DF-EL}
  \forall ~ \bj \in \Z^2,\quad 
  (P_h u)_\bj = \dfrac{\partial L}{\partial y} \left(x_\bj,u_\bj,(\nabla_h u)_\bj\right)
  -
  \left (
    \div_h \mathbf{q}\right)_\bj =0
  \\ \text{with}\quad \mathbf{q\in\mathcal{V}},\quad 
  \mathbf{q}_\bj = 
  \dfrac{\partial L}{\partial z} (x_\bj,u_\bj,(\nabla_h u)_\bj)
  \quad \forall ~\bj\in\Z^2.
\end{align*}
The  differential form for the discretisation of the PDE \eqref{EL} is,
\begin{displaymath}
  \text{find} \quad 
  u \in X_h
  \quad  \text{such that} \quad 
  P_h(u) = 0.
\end{displaymath}
We can define the discrete Lagrangian
$\LL_h:~X_h\longrightarrow \R$  for $u\in X_h$ by,
\begin{equation*}
  % \label{eq:DF-LLh}
  \forall ~ \bj\in\Z^2,\quad 
  \bigl(\LL_h u\bigl)_\bj = \sum_{\bj\in J} L\bigl(x_\bj,u_\bj,(\nabla_h u)_\bj\bigl)~h^2\,.
\end{equation*}
It provides a discrete variational embedding for $\LL$.
The associated discrete variational formulation
of the problem is:
% The associated discrete problem under vraitional formd is, 
\begin{equation*}
  % \label{eq:DF-ELh}
  {\rm find }\, u\in X_h \quad {\rm such~ that} \quad 
  \quad 
  D\LL_h(u)(v)=0
  \quad 
  {\rm for~ any~ } ~v\in X_h.
\end{equation*}
We conserve in this framework the coherence result enunciated in theorem \ref{thm:coherence-df}. 
It is similarly the consequence of  the discrete Green-Gauss 
formula \eqref{eq:DF-green}. Precisely a solution 
to the discrete variational formulation of the problem
satisfies for all $v\in X_h$,
% We conserve in this framework the coherence result enunciated in theorem \ref{thm:coherence-df}, and coherence is induced by the discrete Green-Gauss Formula \eqref{eq:DF-green}. Precisely a solution to the discrete problem under variational form satisfies for all $v_h\in X_h$,
\begin{displaymath}
  \sum_{\bj\in J} \dfrac{\partial  L }{\partial y}
  \left (x_\bj,u_\bj,(\nabla_h u)_\bj \right )v_\bj~h^2
  +
  \sum_{\bj\in J} \dfrac{\partial  L }{\partial z}
  \left (x_\bj,u_\bj,(\nabla_h u)_\bj\right )\cdot (\nabla_h v)_\bj  ~h^2= 0.
\end{displaymath}
With the discrete Green-Gauss formula \eqref{eq:DF-green} we get:
\begin{displaymath}
  \sum_{\bj\in J} \left(
    \dfrac{\partial  L }{\partial y}
    (x_\bj,u_\bj,\nabla_\bj u) 
    -(\div_h\mathbf{q})_\bj
  \right) v_\bj~h^2
  = 0
\end{displaymath}
and we exactly  recover the discrete differential formulation of
the problem.
% and we exactly  recover the discretite problem, differential form.
%% 
%% 
%% 
%% 
%% 
%% 
%% 
%% 
%% 
%% 
%% 
%%%%%%%%%%%%%%%%%%%%%%%%%%%% 
%%%%%%%%%%%%%%%%%%%%%%%%%%%% 
\subsection{Finite Volumes}
%%%%%%%%%%%%%%%%%%%%%%%%%%%% 
%%%%%%%%%%%%%%%%%%%%%%%%%%%% 
%% 
%% 
%% 
%% 
%% 
%% 
%% 
%% 
\label{section:fv}
We focus in this section on the classical finite volume method (as presented \textit{e.g.} in \cite{gallouet}) for the Poisson problem \eqref{eq:pb}.
We consider a mesh $\msh$ of the domain $\Omega$ as in definition \ref{def:msh}. Relatively to this mesh we assume that we can build two sets of points: cell centres $(x_ K)_{ K \in\msh}$  and boundary face centres $(x_ e)_{ e \in\bfc}$ that satisfy:
\begin{align*}
  % \label{eq:fv-msh-adm1}
  \forall ~ K \in\msh, \quad \forall~  e \in\bfc:\quad &x_ K \in K,~x_ e \in e .
  \\
  \notag
  \forall~  e \in\ifc:\quad& e = K_1\vert K_2, ~ [x_{ K_1},x_{ K_2}] \perp  e,~
  \\
  \label{eq:fv-msh-adm3}
  \forall~  e \in\bfc:\quad& e = K \vert \partial\Omega, ~ [x_ e,x_ K ] \perp  e .
\end{align*}
These two conditions are referred to as \textit{admissibility} conditions. They impose a strong constraint on the mesh $\msh$.
Distances $(d_ e)_{ e \in\fc}$ across
the faces are defined as follows:
\begin{align*}
  \notag
  \forall~  e = K_1\vert K_2 \in\ifc:&~d_ e =\vert x_{ K_1}x_{ K_2}\vert,
  \\
  \forall~  e = K \vert\partial\Omega \in\bfc:&~d_ e =\vert x_{ K }x_ e \vert .
\end{align*}
\subsubsection{Discrete differential embedding}
We consider the settings in section \ref{sec:fv-div}: 
$X=[{\rm H}^1(\Omega)]^d$, $Y={\rm L}^2(\Omega)$, $X_h = \R^{\# \fc}$ and 
$Y_h= \R^{\# \msh}$.
We recall that the projections $\pi_1:~X\rightarrow X_h$ and $\pi_2:~Y\rightarrow Y_h$
are the normal component mean values on the mesh faces and the mean values on the mesh cells respectively.
\\
The finite volume divergence $\div_h:~Y_h\rightarrow  X_h$ is defined in equation (\ref{eq:divh-fv}) that we recall,
\begin{displaymath}    
  \div_K \pp_h = \dfrac{1}{|K|} \sum_{e\in\fc, e \subset \partial K} p_e |e| ~\bn_e\cdot \bn{_{K, e }},
\end{displaymath}
with the same notation $\div_K \pp_h := (\div_h \pp_h)_K$.
\\ \\
The flux operator $\mathcal{F}:~{\rm H}^2(\Omega)\rightarrow \R^{\# \fc}$ (thus relatively to the mesh $\msh$) is defined as $\mathcal{F}=\pi_1\circ \nabla $ (it consists in averaging the 
gradient of a function over the edges in their normal direction). 
The discrete flux operator is defined as: 
\begin{displaymath}
  \mathcal{F}_h:~u_h=(u_K)_{K\in\msh}\in \R^{\#\msh} \mapsto (\mathcal{F}_e u_h)_{e\in\fc} \in \R^{\#\fc},
\end{displaymath}
with,
\begin{align*}
  % \label{eq:fv-flx1}
  \forall ~e=K_1\vert K_2\in\ifc~:\quad &
  \mathcal{F}_{e} u_h=\dfrac{u_{K_2}-u_{K_1}}{d_e}\bn_{K_1,e}\cdot\bn_e,
  \\
  \notag
  \forall ~e= K \vert\partial\Omega\in\bfc~:\quad &
  \mathcal{F}_e u_h=-\dfrac{u_{K}}{d_e}\bn_{K,e}\cdot\bn_e.
\end{align*}
Numerical fluxes across edges (and according to their intrinsic orientation) thus are computed  using  a finite difference scheme.
Note that the Dirichlet boundary condition has implicitly being taken into account when defining the numerical fluxes on the boundary faces.
This provides a discrete embedding for the flux operator $\mathcal{F}$:
$$
\xymatrix{
  {\rm H}^2(\Omega)  \ar[r]^{ \mathcal{F}} \ar[d]_{\pi_3}  & \R^{\#\fc}  \\
  \R^{\#\msh} \ar[ur]_{\mathcal{F}_h} & 
}
$$
where the projection $\pi_3$   is defined as $(\pi_3 u)_K=u(x_K)$, {\it i.e.} as the values of the function $u$ at each cell centre $x_K$. 
\\ \\
The discrete Laplace operator $\Delta_h$ %:~\R^{\#\msh} \rightarrow \R^{\#\msh}$
is defined as,
\begin{align*}
  \Delta_h:~\R^{\#\msh} \rightarrow \R^{\#\msh},\quad \Delta_h = \div_h \circ \mathcal{F}_h.
\end{align*}  
For $f\in {\rm L}^2(\Omega)$, the discrete differential embedding of $P(u)=-\Delta u -f$
then is $P_h(u_h) = -\Delta_h u_h -\pi_2 f$ for $u_h\in \R^{\#\msh}$.
The differential formulation for the discrete Poisson problem is, 
\begin{equation}
  \label{eq:fv-scheme}
  {\rm find }~u_h\in \R^{\#\msh} \quad {\rm so~ that}\quad 
  P_h(u_h) = -\Delta_h u_h - \pi_2 f=0.
\end{equation}
Moreover we have the following \textit{discrete Green-Gauss formula}:
for all $p=(p_e)_{e\in\fc}\in \R^{\#\fc}$ and for all $u_h=(u_K)_{K\in\msh}\in \R^{\#\msh}$,
\begin{equation}
  \label{eq:fv-green}
  \sum_{ K \in\msh}
  \left(
    \div_K p \right) u_K |K|
  = -\sum_{e\in\fc}
  p_e \mathcal{F}_e u_h |e|d_e.
\end{equation}
\subsubsection{Discrete variational embedding, coherence}
In the continuous case, the diffusion energy $\int_\Omega\vert \nabla u\vert ^2/2\,\dx$ is part
of the Lagrangian functional $\LL$.
In the framework of finite volume method, no proper \textit{discrete gradient} is available, but only numerical fluxes in the normal direction to the mesh faces. 
Thus only the normal component (and not the tangential one)
of some  \textit{discrete gradient} on the mesh faces is approximated. 
\\
  We recall that $\pi_2 f = (f_K)_{K\in \msh}$ with $f_K=\int_K f\, \dx / |K|$, see  section \ref{sec:fv-div}.
  We introduce the discrete Lagrangian functional $\LL_h :~\R^{\#\msh}\rightarrow\R$ as
\begin{displaymath}
  \LL_h(u_h)=
  \dfrac{1}{2}\sum_{ e \in\fc} (\mathcal{F}_{e} u_h)^2 ~|e|d_e
  ~-~
  \sum_{ K \in\msh} f_K u_K~\meas{K}.
\end{displaymath}
The functional $\LL_h$ defines a discrete variational embedding of $\LL$.
The variational form for the finite volume discrete Poisson problem is,
\begin{equation}\label{eq:fv-var}
  {\rm find}\quad u_h\in \R^{\#\msh} \quad {\rm such~ that} \quad 
  \forall ~v_h\in \R^{\#\msh},\quad   D\LL_h(u_h)(v_h)=0.
\end{equation}
\begin{theorem}%[Coherence]
  The discrete variational and differential embeddings of the Poisson problem using  the finite volume method  are coherent. 
  Precisely, the discrete differential formulation \eqref{eq:fv-scheme} and
  discrete variational formulation \eqref{eq:fv-var} for the Poisson problem have the same solutions.
\end{theorem}
\begin{proof} 
  Let $u_h$ satisfy \eqref{eq:fv-var}, we have by differentiating $\LL_h$:
  for all $u_h,v_h \in \R^{\#\msh}$,
  \begin{align*}
    D\LL_h(u_h)(v_h)=
    \sum_{e\in\fc}      \bigl(\mathcal{F}_e u_h\bigl) \bigl(\mathcal{F}_e v_h\bigl) |e|d_e
    ~-~
    \sum_{ K \in\msh} (\pi_2 f)_K v_K\meas{K}=0.
  \end{align*}
  Using the discrete Green-Gauss formula \eqref{eq:fv-green}, we  get for all $u_h,v_h \in\R^{\#\msh}$,:
  \begin{displaymath}
    - \sum_{ K \in\msh} \left(\div_K (\mathcal{F}_h u_h)\right)v_K |K|
    - \sum_{ K \in\msh} (\pi_2f)_K v_K |K|
    =0,
  \end{displaymath}
  which is equivalent with \eqref{eq:fv-scheme}.
\end{proof}
%% 
%% 
%% 
%% 
%% 
%% 
%% 
%%%%%%%%%%%%%%%%%%%%%%%%%%%%%%%%%%%%%%%%%%%%%%%% 
%%%%%%%%%%%%%%%%%%%%%%%%%%%%%%%%%%%%%%%%%%%%%%%% 
\section{Hamiltonian calculus of variations and mixed formulations}
%%%%%%%%%%%%%%%%%%%%%%%%%%%%%%%%%%%%%%%%%%%%%%%% 
%%%%%%%%%%%%%%%%%%%%%%%%%%%%%%%%%%%%%%%%%%%%%%%% 
%% 
%% 
\label{part3}
In this section let  $ L$ % :~\Omega\times \R\times \R^d$
be an 
admissible Lagrangian function as defined in section \ref{sec:cl-euler-lag}. 
We recall the link between Hamiltonian and Lagrangian systems
in section \ref{sec:ham}. 
We will stress here  the relationships between mixed formulations and
discrete embedding of Hamiltonian systems in section
% \todo{cette section ne fait plus ca }\ref{sec:disc-ham} and 
\ref{sec:mfd}.
%% 
%% 
%% 
%% 
%% 
%% 
%%%%%%%%%%%%%%%%% new Charles 17 juillet 2012 %%%%%%%%%%%%%%%%%%%%%5
%% 
%% 
%% 
%% 
%% 
%% 
%% 
%%%%%%%%%%%%%%%%%%%%%%%%%%%% 
%%%%%%%%%%%%%%%%%%%%%%%%%%%% 
\subsection{Hamiltonian formulation}
%%%%%%%%%%%%%%%%%%%%%%%%%%%% 
%%%%%%%%%%%%%%%%%%%%%%%%%%%% 
%% 
%% 
\label{sec:ham}
%% 
%% 
%% 
%% 
%% 
%% 
%% 
%% 
%% 
%% 
%% 
%\subsubsection{Hamiltonian formulation}
%% 
%% 
%% 
%% 
\label{sub:ham}
\begin{definition}[Legendre property]
  We say that  $ L $ satisfies the Legendre property if the mapping 
  $z\mapsto \dfrac{\partial  L }{\partial z}(x,y,z)$ is a bijection on $\R^d$ for any $x\in \Omega$ and any $y\in \R$. 
\end{definition}
If $ L $ satisfies the Legendre property, the following function  $g:~\Omega\times \R\times \R^d\rightarrow  \R^d$ is well defined: 
\begin{displaymath}
  z=g(x,y,\pp)\quad {\rm with } \quad  \pp=\di {\partial  L \over \partial z}(x,y,z).
\end{displaymath}
Let us consider $\pp=\di {\partial  L \over \partial z}(x,y,z)$ 
as a new variable, 
then,
\begin{displaymath}
  \pp= \di {\partial  L \over \partial z}\bigl(x,y,g(x,y,\pp)\bigl)
  \quad  {\rm  and} \quad
  g\bigl(x,y, \di {\partial  L \over \partial z}(x,y,z)\bigl) = z\,.
\end{displaymath}
\begin{definition}[Hamiltonian]
  \label{def:hamilton}
  Let $ L $ satisfy the Legendre property.
  The Hamiltonian 
  $ H :\Omega\times\R \times \R^d \rightarrow \R$ 
  associated to $ L $ is:
  \begin{equation*}
    H (x,y,\pp)=\pp \cdot g(x,y,\pp)- L (x,y,g(x,y,\pp)).
  \end{equation*}
  We introduce two different definitions for the Hamiltonian functional $\HH:~{\rm Dom}(\HH)\subset  {\rm L}^2(\Omega)\times \left[ {\rm L}^2(\Omega)\right]^d\rightarrow \R$ associated to $ H $:
  % with domain ${\rm Dom}(\HH)\subset  {\rm L}^2(\Omega)\times \left[ {\rm L}^2(\Omega)\right]^d$:
  \begin{itemize}
  \item [1.] Primal Hamiltonian, ${\rm Dom}(\HH)={\rm H}^1(\Omega)\times \left[ {\rm L}^2(\Omega)\right]^d$,
    \begin{equation}
      \label{eq:HH-1}
      \HH(u,\pp)=\int_{\Omega} \pp\cdot \nabla u  - H (x,u,\pp)\,\dx.
    \end{equation}
  \item [2.] Dual Hamiltonian, ${\rm Dom}(\HH)={\rm L}^2(\Omega)\times {\rm H}_{\div}(\Omega)$, 
    \begin{equation}
      \label{eq:HH-2}
      \HH(u,\pp)=\int_{\Omega} -\div(\pp) u  - H (x,u,\pp)\,\dx.
    \end{equation}
  \end{itemize}
\end{definition}
\begin{proposition}
  % We consider a space of variations $V\times W\subset {\rm Dom}(\HH)$. 
  According to definition \ref{def:diff}, the Hamiltonian functional $\HH$ is differentiable at point $(u,\pp)\in {\rm Dom}(\HH)$ if 
  \begin{equation*}
    \dfrac{\partial H }{\partial y} (x,u,\pp)\in {\rm L}^2(\Omega)
    \quad 
    {\rm 
      and }
    \quad
    \dfrac{\partial H }{\partial \pp } (x,u,\pp)\in\left[ {\rm L}^2(\Omega)\right]^d .
  \end{equation*} 
  In such a case we have, for $(v,\qq)\in {\rm Dom}(\HH)$:
  \begin{itemize}
  \item In the primal case: 
    \begin{equation}
      \label{eq:der-HH}
      D\HH(u,\pp)\cdot(v,\qq)=
      \int_\Omega \left[
        \qq\cdot\left(
          \nabla u - \dfrac{\partial H}{\partial \pp}(x,u,\pp)
        \right)
        +\nabla v\cdot \pp
        - v \dfrac{\partial H}{\partial y}(x,u,\pp)
      \right]\dx.
    \end{equation}
  \item In the dual case:
    \begin{displaymath}
      D\HH(u,\pp)\cdot(v,\qq)=
      \int_\Omega \left[
        -\div(\qq)u
        -\qq\cdot\dfrac{\partial H}{\partial \pp}(x,u,\pp)   
        -v \left(\div\pp
          +\dfrac{\partial H}{\partial y}(x,u,\pp)
        \right)
      \right]\dx.
    \end{displaymath}
  \end{itemize}
\end{proposition}
\begin{definition}[Extremals]
  \label{def:ext-hh}
  Let us consider a space of variation $V\times W\subset {\rm Dom}(\HH)$.
  We say that $(u,\pp)\in {\rm Dom}(\HH)$ is an extremal for $\HH$  relatively to $V\times W$ if $\HH$ is differentiable at point $(u,\pp)$ and:
  \begin{displaymath}
    \forall~(v,\qq)\in V\times W,\quad D \HH(u,\pp)\cdot(v,\qq) = 0.
  \end{displaymath}
\end{definition}
\begin{theorem}[Hamilton's least action principle] 
  \label{thm:hamilton-lap}
  Let $(u,\pp)\in {\rm Dom}(\HH)$ be an extremal for $\HH$ relatively to $V\times  W$. 
  Assume moreover that:
  \begin{itemize}
  \item in the primal case: $\pp\in {\rm H}_{\div}(\Omega)$,
    $V_0=\{v\in V, v=0 {\rm  ~on~ } \partial\Omega\}$ 
    is dense in ${\rm L}^2(\Omega)$ 
    and $W$ is dense in $\left[{\rm L}^2(\Omega)\right]^d$,
  \item in the dual case:  $u\in {\rm H}^1(\Omega)$,
    $V$ is dense in ${\rm L}^2(\Omega)$ and 
    $W_0=\{\qq\in W, \qq\cdot n=0 {\rm  ~on~ } \partial\Omega\}$ 
    is dense in $\left[{\rm L}^2(\Omega)\right]^d$.
  \end{itemize}
  Then $(u,\pp)$ is a solution of the  {\sl  Hamiltonian system}:
  \begin{equation}
    \label{eq:sysham}
    \left \{
      \begin{array}{ll}
        {\di \div \pp} & {\di = -\dfrac{\partial  H }{\partial y }(x,u,\pp) }
        \\[8pt]
        {\di \nabla u}  & {\di =\dfrac{\partial  H }{\partial \pp}(x,u,\pp)} .
      \end{array}
    \right .
  \end{equation}
\end{theorem}
\begin{proof}
  Let us consider the case of the primal definition of the Hamiltonian functional $\HH$.
  Since $\pp\in {\rm H}_{\div}(\Omega)$, using the Green formula in \eqref{eq:der-HH} gives: $\forall~ (v,\qq) \in V\times W$,
  \begin{displaymath}
    \int_{\Omega}  \Bigl(
    -\bigl(\div \pp + \frac{\partial H}{\partial y}(x,u,\pp)\bigl)\,v
    +\qq \cdot\bigl(\nabla u-\frac{\partial H}{\partial \pp}(x,u,\pp) \bigl)
    \Bigl)\dx ~+~\int_{\partial\Omega} v~ \pp\cdot\bn~ds  = 0.
  \end{displaymath}
  The boundary integral vanishes for $v\in V_0$.
  We recover \eqref{eq:sysham} by density of $V_0$ in ${\rm L}^2(\Omega)$ and of $W$ in $\left[{\rm L}^2(\Omega)\right]^d$. 
\end{proof}
\begin{corollary}[Lagrangian and Hamiltonian formulations]
  The solutions $(u,\pp)$ of the Hamiltonian system \eqref{eq:sysham} are exactly the solutions of the Euler-Lagrange equation \eqref{EL} under the condition
  \begin{equation*}
    \pp = \di {\partial L\over \partial z}(x,u, \nabla u).
  \end{equation*}
\end{corollary}
\subsubsection*{Application to the Poisson problem}
We consider the
Poisson problem (\ref{eq:pb}). 
We recall that the 
Lagrangian function associated with this problem is 
$$ L (x,y,z)=\frac{1}{2} z\cdot z-f(x)y.$$
The Legendre property is clearly satisfied by $ L $.
We introduce the new variable $\pp=z$ and
the function $g$ is given by $g(x,y,\pp)= \pp$.
A Hamiltonian for the Poisson problem is then given by,
\begin{equation}
  \label{def:H-poisson}
  H (x,y,\pp)=\pp\cdot\pp -  L (x,y,g(x,y,\pp)) =
  \dfrac{1}{2}\pp \cdot \pp + f(x)y.
\end{equation}
The Hamiltonian system \eqref{eq:sysham} 
associated with \eqref{def:H-poisson} 
is the mixed formulation of the Poisson problem (\ref{eq:pb}):
\begin{equation}
  \label{mixte}
  \left\{
    \begin{array}{ll}
      {\di -\div \pp} &= {\di f}
      \\[6pt]
      {\di \nabla u}&={\di \pp}\,.
    \end{array}
  \right.
\end{equation}
Applying theorem \ref{thm:hamilton-lap}, one obtains that the weak solutions of the Poisson problem in its mixed form
\eqref{mixte}
exactly are extremals for the Hamiltonian functional $\HH$ in \eqref{def:H-poisson}. Precisely:
\begin{itemize}
\item Primal form \eqref{eq:HH-1} of $\HH$. 
  Consider an extremal $(u,\pp)\in {\rm H}_0^1(\Omega)\times \left[{\rm L}^2(\Omega)\right]^d$ of $\HH$ relatively to the space of variations  $V\times  W={\rm H}_0^1(\Omega)\times \left[{\rm L}^2(\Omega)\right]^d$. 
  An extremal exactly is a  solution for the primal weak formulation of the mixed Poisson equation: find $(u,\pp)\; \in {\rm H}_0^1(\Omega)\times \left[{\rm L}^2(\Omega)\right]^d$ such that,
  \begin{displaymath}
    % \label{eq:T.1.2.2.2.0}
    \left\{
      \begin{array}{lll}
        -\int_\Omega \pp\cdot\nabla v ~\dx =-\int_\Omega fv~\dx \quad &\quad \forall~ v \in {\rm H}_0^1(\Omega)
        \\[5pt]
        \int_\Omega (\pp-\nabla u)\cdot\qq ~\dx  =0 \quad &\quad \forall~ \qq \in \left[{\rm L}^2(\Omega)\right]^d\,.
      \end{array}
    \right. 
  \end{displaymath}\\
\item  Dual form \eqref{eq:HH-2} of $\HH$. 
  Consider an extremal $(u,\pp)\in {\rm L}^2(\Omega)\times {\rm H}_{\div}(\Omega)$  of $\HH$ relatively to the space of variations  $V\times  W = {\rm L}^2(\Omega)\times {\rm H}_{\div}(\Omega)$. 
  An extremal exactly is a  solution for the dual weak formulation of the mixed Poisson equation that reads:  find $(u,\pp)\; \in {\rm L}^2(\Omega)\times {\rm H}_{\div}(\Omega)$ such that
  \begin{equation*}
    % \label{eq:T.1.2.2.2}
    \left\{
      \begin{array}{lll}
        \int_\Omega(\div \pp+f)v~\dx=0 \quad&\quad  \forall~ v \in {\rm L}^2(\Omega)
        \\[5pt] 
        \int_\Omega\pp\cdot\qq~\dx+\int_\Omega u\div \qq~\dx=0 \quad&\quad \forall~ \qq \in {\rm H}_{\div}(\Omega)\,.
      \end{array}
    \right. 
  \end{equation*}
\end{itemize}
\subsection{Coherence}
\label{sec:disc-ham}
The definition of the discrete  differential embedding
in section \ref{part:embedding}
applies to the Hamiltonian system where $P$ is given by
\begin{equation*}
  P(u,\pp)=
  \left(
    \begin{array}{ll}
      {\di \div \pp+\dfrac{\partial  H }{\partial y }(x,u,\pp) }
      \\[8pt]
      {\di \nabla u -\dfrac{\partial  H }{\partial \pp}(x,u,\pp)} 
    \end{array}
  \right) .
\end{equation*}
The definition \ref{def:disc-var-emb} of the discrete
variational embedding also applies to Hamiltonian system by replacing
$\LL$ by $\HH$.
The definition of coherence for the discretisation of Hamiltonian systems is
the same as definition \ref{def:coherence} for Lagrangian systems.
% in definition \ref{def:coherence}.
%% 
\begin{definition}
  \label{def:coherence-H}
  Let us consider a  discrete differential embedding  of the mixed problem 
  \eqref{eq:sysham}. If the discrete problem solutions exactly are extremals of 
  a discrete Hamiltonian functional $\HH_h$ that moreover also is a 
  discrete variational embedding of $\HH$, then we say that we have coherence.
  \\
  In case of coherence we then have the following commutative diagram:
  \begin{eqnarray*}
    \begin{CD}
      {{(u,\pp) \mapsto \HH(u,\pp)}} &   @>{{\rm disc.~ var.~ emb.}}>> & 
      { {(u_h,\pp_h) \mapsto \HH_h(u_h,\pp_h)}}\\
      @V{{\rm L.A.P.}}VV & \hfill & @VV{{\rm  disc.~ L.A.P. }}V\\
      (u,\pp) \quad{\rm  solution~ of~ PDE~ \eqref{eq:sysham}}  & @>{{\rm disc.~ diff.~ emb.}}>> &   
      (u_h,\pp_h) \quad {\rm solution~ of~ the~ discrete~ PDE}\\
      {\rm Hamiltonian~ system}  & @.& {{\rm  disc.~ Hamiltonian~ system}}
    \end{CD}
  \end{eqnarray*}
  where L.A.P. stands for least action principle.
\end{definition}
\begin{remark}
  In section \ref{sec:coherence}
  it was shown that the coherence for conforming finite element
  naturally derives from the method definition.
  The same conclusion also holds 
  for conforming mixed finite elements.
  The discrete Hamiltonian in that case is the restriction of the Hamiltonian $\HH$ 
  to the finite element space.
\end{remark}
%% 
%% 
%% 
%% 
%% 
%%%%%%%%%%%%%%%%%%%%%%%%%%%% 
%%%%%%%%%%%%%%%%%%%%%%%%%%%% 
\subsection{Mimetic Finite Differences}
%%%%%%%%%%%%%%%%%%%%%%%%%%%% 
%%%%%%%%%%%%%%%%%%%%%%%%%%%% 
%% 
%% 
%% 
%% 
%% 
%% 
%% 
\label{sec:mfd}
We consider the mixed formulation  \eqref{mixte} of the Poisson problem together with a homogeneous Dirichlet condition $u=0$ on $\partial\Omega$.
% The dual weak formulation \eqref{eq:T.1.2.2.2} is adopted. 
The scalar products on ${\rm L}^2(\Omega)$ and on $\left[{\rm L}^2(\Omega)\right]^d$ are respectively denoted, for $u$, $v\in {\rm L}^2(\Omega)$ and for $\pp$, $\qq \in \left[{\rm L}^2(\Omega)\right]^d$,
\begin{displaymath}
  \left(u,v\right) = \int_\Omega uv~\dx,\quad \left[\pp,\qq\right] =
  \int_\Omega \pp\cdot \qq ~\dx.
\end{displaymath}
The Green-Gauss formula rewrites as, for all $u\in {\rm H}^1_0(\Omega)$
and all $\pp\in {\rm H}_{\div}(\Omega)$,
\begin{displaymath}
  \left[\pp,\nabla u\right] = -\left(\div \pp, u \right).
\end{displaymath}
\\
In the Mimetic Finite Differences (MFD) framework, a discrete flux operator $\mathcal{F}_h~:Y_h\rightarrow X_h$ is defined as (minus) the adjoint of the finite volume discrete divergence (see section \ref{sec:fv-div}) after the introduction of a scalar product on $X_h$  that is consistent with $[\cdot,\cdot]$.
We refer to \cite{brezzi-2005} for the MFD discretisation of diffusion problems.
\subsubsection{Discrete differential embedding}

A mesh $\msh$ of the domain $\Omega$ is considered as in definition \ref{def:msh}.
The space $P^0(\msh)$ of the piecewise constant functions over the mesh cells is considered and identified with $\R^{\#\msh}$.
Since $P^0(\msh)\subset {\rm L}^2(\Omega)$, the ${\rm L}^2$ scalar product on $P^0(\msh)$ is available.% , 
\\
The notations in section \ref{sec:fv-div} for the finite volume divergence 
are considered:   
$X_h=\R^{\#\fc}$ and
$\pi_1:~[{\rm H}^1(\Omega)]^d\rightarrow X_h$,
$\pi_2:~{\rm L}^2(\Omega)\rightarrow P^0(\msh)$ are the projections 
in the diagram (\ref{eq:disc-div-fv}). 
We adopt the following alternative (but equivalent) definition for the finite volume divergence
in the diagram (\ref{eq:disc-div-fv}). 
  We introduce the space $\tilde{X}_h$:
  \begin{eqnarray*}
    \tilde{X}_h = 
    \left\{
      p_{K,e}\quad \text{for} \quad K\in \msh 
      \quad \text{and for} \quad e\in\fc
      \quad \text{so that} \quad e\subset \partial K
    \right.
    \\
    \left.
      \text{that satisfy}\quad
      p_{K_1,e}+p_{K_2,e}=0 \quad {\rm  if }\quad e=K_1|K_2 
    \right\}.
  \end{eqnarray*}
  An element $\mathbf{p}\in\tilde{X}_h$ is given by one numerical flux on each external face and by two opposite numerical fluxes per internal face.
  Obviously,  $\tilde{X}_h$ is isomorphic to $X_h$.
With this identification we get the new commutative diagram for the discrete divergence,
\begin{equation*}
  \begin{CD}
    \left[{\rm H}^1(\Omega)\right]^d &   @>{\div}>> & {\rm L}^2(\Omega)
    \\
    @V{\tilde{\pi}_1}VV & \hfill & @VV{\pi_2}V
    \\
    \tilde{X}_h &   @>{\div_h}>> &  P^0(\msh) 
  \end{CD}    
\end{equation*}
where $\tilde{\pi}_1$ is given by $\tilde{\pi}_1 \pp=(p_{K,e})$ with $p_{K,e}=\int_e \pp\cdot \bn_{K,e} \dl / |e|$ the mean flux of $\pp$ across the face $e$ according to the unit normal to $e$ pointing outwards $K$. The discrete divergence within this framework has the following expression (to be compared to \eqref{eq:divh-fv}), $\div_h:~ \pp =(p_{K,e}) \in  \tilde{X}_h \mapsto (\div_K \pp)_{K\in \msh} \in P^0(\msh)$:
$$
\div_K \pp = \dfrac{1}{|K|} \sum_{e\in\fc, e \subset \partial K} p_{K,e} |e| .
$$
The definition of a scalar product on $\tilde{X}_h$ is not obvious. Let us
consider $ K \in\msh$ and denote $\tilde{X}_h^ K $ the restriction of $\tilde{X}_h$ to
$ K $. We suppose that a cell scalar product $[\cdot,\cdot]_K $ is
given on each $\tilde{X}_h^ K \in\msh$ and that the scalar product on $\tilde{X}_h$ decomposes as:
\begin{equation}
  \label{eq:mfd-scal-flx}
  \forall~ \pp_h,\qq_h\in \tilde{X}_h:\quad 
  [\pp_h,\qq_h]_h = \sum_{ K \in\msh} [\pp_h,\qq_h]_{ K },
\end{equation}
A way to define the elemental scalar product \eqref{eq:mfd-scal-flx}
is to introduce a lifting operator $\mathcal{R}_K: \tilde{X}_h^ K \longrightarrow \left[{\rm L}^2(K)\right]^d$
and then to define:
\begin{equation}
  \label{def[]}
  [\pp_h,\qq_h]_ K  = \int_ K  \mathcal{R}_K(\pp_h)\cdot\mathcal{R}_K(\qq_h) \dx.
\end{equation}

For more details on the construction of $\mathcal{R}_K$,
we refer to \cite{brezzi-2005}. The present definitions are sufficient for our purpose.
Relatively to the scalar product
% \eqref{eq:mfd-prod-scal} and
\eqref{eq:mfd-scal-flx}, the discrete flux
operator $\mathcal{F}_h: P^0(\msh) \longrightarrow \tilde{X}_h$ is defined as  (minus) the adjoint of the
discrete divergence: $\mathcal{F}_h=-\div_h^\star$.
It is uniquely determined by,
\begin{displaymath}
  \forall ~u_h,~\pp_h\in P^0(\msh)\times \tilde{X_h}:\quad 
  \left [\pp_h, \mathcal{F}_h u_h \right ]_{h}
  = - (\div_h \pp_h, u_h).
\end{displaymath}
The discrete differential embedding for the mixed Poisson problem (\ref{mixte}) using the MFD method then is defined by $P_h:~P^0(\msh)\times  \tilde{X}_h \rightarrow P^0(\msh)\times  \tilde{X}_h$:
\begin{displaymath}
  P_h(u_h,\pp_h)=  
  \left(
    \begin{array}{ll}
      {\di -\div_h \pp_h   - \pi_2 f}
      \\[5pt]
      {\di \pp_h-\mathcal{F}_h u_h  }
    \end{array}
  \right) . 
\end{displaymath}
The discretisation of the  mixed Poisson problem \eqref{mixte} is
% (differential form)
:
find $u_h\in P^0(\msh) $ and $\pp_h\in \tilde{X}_h$ such that,
\begin{equation}
  \label{eq:mfd-disc-prob}
  P_h(u_h,\pp_h)=  0.
\end{equation} 
\subsubsection{Discrete variational embedding, coherence}
The Hamiltonian $H$ for the Poisson problem is given in equation  \eqref{def:H-poisson}.
The associated  Hamiltonian functional $\HH$ 
with the primal definition \eqref{eq:HH-1},
${\rm Dom}(\HH)={\rm H}^1(\Omega)\times \left[ {\rm L}^2(\Omega)\right]^d$,
% on ${\rm L}^2(\Omega)\times {\rm H}_{\div}(\Omega)$ is:
\begin{align*}
  \HH(u,\pp) 
  &= \int_\Omega  \pp\cdot\nabla u \,\dx
  - \dfrac{1}{2}\int_\Omega  \pp\cdot\pp\,\dx
  - \int_\Omega fu\,\dx
  \\
  & = \left[\pp,\nabla u\right]
  - \dfrac{1}{2}\left[\pp,\pp\right]
  - \left(u,f\right).
\end{align*}
We therefore define the discrete Hamiltonian $\HH_h:~ P^0(\msh) \times \tilde{X}_h\rightarrow \R$ as,
\begin{displaymath}
  % \label{eq:mfd-hamilton-disc}
  \HH_h(u_h,\pp_h)=\left[\pp_h,\mathcal{F}_h u_h\right]_h
  - \dfrac{1}{2}\left[\pp_h,\pp_h\right]_h
  - \left(u_h, \pi_2 f\right)_h.
\end{displaymath}
It provides the following discrete variational embedding,
$$
\xymatrix{
  {\rm H}^1_0(\Omega)\times \left[{\rm H}^1(\Omega)\right]^d  \ar[r]^{\quad \qquad\HH} \ar[d]_{\pi_2\times \tilde{\pi}_1}  & \R  \\
  P^0(\msh) \times \tilde{X}_h  \ar[ur]_{\quad \HH_h} & 
}
$$
The variational form for the MFD discrete mixed Poisson problem is:
find $(u_h,\pp_h)\in P^0(\msh) \times \tilde{X}_h$ such that,
\begin{equation}\label{eq:mfd-var}
  \forall ~(v_h,\qq_h)\in P^0(\msh) \times \tilde{X}_h ,\quad   
  D\HH_h(u_h,\pp_h)(v_h,\qq_h)=0.
\end{equation}
\begin{theorem}
  The MFD discrete differential formulation 
  \eqref{eq:mfd-disc-prob} and variational formulation 
  \eqref{eq:mfd-var} for the mixed Poisson problem are equivalent.
  Then the MFD discretisation for the mixed Poisson problem
  is coherent.
\end{theorem}
\begin{proof}
  Differentiating $\HH_h$ gives:
  \begin{displaymath}
    D\HH_h(u_h,\pp_h)(v_h,\qq_h) =    
    \left[\pp_h,\mathcal{F}_h v_h\right]_h
    + \left[\mathcal{F}_h u_h,\qq_h\right]_h
    - \left[\pp_h,\qq_h\right]_h
    - \left(\pi_2 f,v_h\right)_h
    .
  \end{displaymath}
  Using that %$\mathcal{F}_h=-\mathcal{F}_h^\star$ 
  $\mathcal{F}_h=-\div_h^\star$ relatively to the scalar product $[\cdot,\cdot]_h$
  we obtain,
  \begin{displaymath}
    D\HH_h(u_h,\pp_h)(v_h,\qq_h) =    
    \left(-\div_h \pp_h-\pi_2 f, v_h\right)_h
    + \left[\mathcal{F}_h u_h-\pp_h,\qq_h\right]_h
    .
  \end{displaymath}
  Therefore singular points for $\HH_h$ exactly are the solutions to equation \eqref{eq:mfd-disc-prob}.
\end{proof}

\section{Conclusion}
In the present paper we studied the properties of the discretisation of PDEs deriving from a variational principle, either Lagrangian or Hamiltonian.
We addressed the following questions. Does the discrete problem also satisfy a variational principle ? If it does, what is the relationship between that variational principle and the one that rules the PDE ?
These questions are analysed by introducing the concepts of discrete variational and discrete differential embeddings and of coherence between these two
% sorts
types of embeddings.
\\
For the Poisson problem, considering several classical methods,
we showed that the discrete Poisson equation is associated to a variational %property.
embedding.
% By considering the Poisson problem as a test case, we showed for several classical methods that the discrete Poisson equation satisfy a variational %property.
% embedding.
A crucial property ensuring coherence for the discrete problems is the following.
The %order 2 
Euler Lagrange PDE involves two  differential operators of order one:
a gradient and a divergence.
The differential embeddings of these two operators must satisfy some duality property.
That property is a discrete analogous of the Green-Gauss formula.
%% 
%% 
%% 
%% 
%% 
%% 
%% 
%% 
%% 
%%%%%%%%%%%%%%%%%%%%%%%%%%%% 
%%%%%%%%%%%%%%%%%%%%%%%%%%%% 
%% 
%% 
%% 
%% 
%% 
%% 
%%%%%%%%%%%%%%%%%%%%%%%%%%%%%%%%%%%%%%%%%% 
% BIBLIO
%%%%%%%%%%%%%%%%%%%%%%%%%%%%%%%%%%%%%%%% 
\bibliographystyle{abbrv}
\bibliography{biblio}

\end{document}